\documentclass[a4paper]{siamltex}

\usepackage{amsmath, amssymb, amsfonts, euscript} 

\newtheorem{remark}[theorem]{{\it Remark\/}}

\newlength{\textequation} 
\DeclareMathOperator{\Cos}{Cos}
\DeclareMathOperator{\Sin}{Sin}
\DeclareMathOperator{\Real}{Re}

\renewcommand{\Re}{\Real}

\newcommand{\C}{\mathbb{C}}
\newcommand{\R}{\mathbb{R}}
\newcommand{\Rp}{\mathbb{R_{+}}}
\newcommand{\N}{\mathbb{N}}

\providecommand{\abs}[1]{\lvert #1 \rvert}
\providecommand{\scal}[2]{\langle #1,#2 \rangle}
\providecommand{\norm}[1]{\lVert #1 \rVert}

\newcommand{\res}[1]{_{\rceil #1}}
\newcommand{\adh}[1]{\, \overline{#1}}

\renewcommand{\geq}{\geqslant}
\renewcommand{\leq}{\leqslant}

\def\vv{\underline{v}}

\def\zz{\underline{z}}

\def\kk{\underline{k}}

\def\loc{_{\text{loc}}}

\def\Tmin{T_{\text{min}}} 

\def\A{\mathcal{A}}

\def\smA{\sqrt{\!  -A}}

\providecommand{\scd}[2]{\scal{#1}{#2}}

\def\dx{\dot{x}}

\def\dz{\dot{z}}
\def\ddz{\ddot{z}}

\def\ddzz{\ddot{\zz}}
\def\dza{\dot{\zeta}}

\def\L{\mathcal{L}} 
\def\D{D} 

\begin{document}


\title{The control transmutation method\\ 
and the cost of fast controls\thanks{This version: February 4, 2004.}
}
\author{Luc Miller\thanks{Centre de Math{\'e}matiques, U.M.R. 7640 du C.N.R.S.,
Ecole Polytechnique, 91128 Palaiseau Cedex, France and 
{\'E}quipe Modal'X, E.A. 3454, B{\^a}timent G, 
Universit{\'e} de Paris X - Nanterre, 
200 Avenue de la R{\'e}publique, 92001 Nanterre Cedex, France 
(miller@math.polytechnique.fr).
This work was partially supported by the ACI grant 
``{\'E}quation des ondes : oscillations, dispersion et contr{\^o}le''.}}

\maketitle

\begin{abstract} 
In this paper, the null controllability in any positive time $T$
of the first-order equation 
(1) $\dx(t)=e^{i\theta}Ax (t)+Bu(t)$ 
($\abs{\theta}< \pi/2$ fixed) 
is deduced from the null controllability in some positive time $L$
of the second-order equation 
(2) $\ddz(t)=Az(t)+Bv(t)$.
The differential equations (1) and (2) are set in a Banach space, 
$B$ is an admissible unbounded control operator,
and $A$ is a generator of cosine operator function.

The control transmutation method 
explicits the input function $u$ of (1)
in terms of the input function $v$ of (2):
$u(t,x)=\int_{\mathbb{R}} k(t,s)v(s)\, ds $,
where the compactly supported kernel $k$ depends on $T$ and $L$ only.
It proves that 
the norm of a $u$ steering the system (1) from an initial state $x_{0}$
to zero grows at most like $\norm{x_{0}}\exp(\alpha_{*} L^{2}/T)$
as the control time $T$ tends to zero.
(The rate $\alpha_{*}$ is characterized independently 
by a one-dimensional controllability problem.)

In the applications to the cost of fast controls 
for the heat equation,
$L$ is the length of the longest ray of geometric optics
which does not intersect the control region.
\end{abstract}







\pagestyle{myheadings}
\thispagestyle{plain}
\markboth{LUC MILLER}{THE CONTROL TRANSMUTATION METHOD} 


\section{Introduction}
\label{sec:intro}

This paper concerns the relationship between the null-con\-trollability 
of the following first and second order controllable systems:
\begin{align}\label{eqfsystcon}
&\dx(t)=e^{i\theta}Ax(t)+Bu(t) & (t\in \Rp ), 
& \quad  x(0)=x_{0}, 
\\ \label{eqsystcon}
&\ddz(t)=Az(t)+Bv(t) & (t\in \R ), 
& \quad z(0)=z_{0},\ \dz(0)=0, 
\end{align}
where 
$x$ and $z$ are the systems trajectories in the Banach space $X$,
$x_{0}$ and $z_{0}$ are initial states,
$u$ and $v$ are input functions with values in the Banach space $U$,
$A$ is an unbounded generator, 
$B$ is an unbounded control operator,
$\theta$ is a given angle in $]-\pi/2,\pi/2[$,
each dot denotes a derivative with respect to the time $t$ 
and $\Rp=[0,\infty)$
(the detailed setting is given in \S\ref{sec:setting}).

Equation (\ref{eqfsystcon}) with $u=0$
describes an irreversible system (always smoothing) 
and we think of it as a parabolic distributed system 
with infinite propagation speed.
Equation (\ref{eqsystcon}) with $v=0$ 
describes a reversible system (e.g. conservative) 
and we think of it as a hyperbolic distributed system 
with finite propagation speed.
For example, if 
$A$ is the negative Laplacian on a Euclidean region and 
the input function is a locally distributed boundary value set by $B$,
then (\ref{eqsystcon}) is a boundary controlled scalar wave equation
and (\ref{eqfsystcon}) with $\theta=0$ is a boundary controlled heat equation 
(\S\ref{sec:appli} elaborates on this example).

This paper presents the {\em control transmutation method} 
(cf.~\cite{Her75} for a survey on transmutations in other contexts)
which can be seen as a shortcut to 
Rusell's famous harmonic analysis method in \cite{Rus73}.
It consists in explicitly constructing controls $u$
in any time $T$ for the heat-like equation (\ref{eqfsystcon})
in terms of controls $v$ in time $L$
for the corresponding wave-like equation (\ref{eqsystcon}),
i.e. $u(t)=\int_{\mathbb{R}} k(t,s)v(s)\, ds $,
where the compactly supported kernel $k$ depends on $T$ and $L$.
It proves that the exact controllability of (\ref{eqsystcon}) in some time $L$
implies the null controllability of (\ref{eqfsystcon}) in any time
with a relevant 
upper bound on the cost of fast controls for (\ref{eqfsystcon}).
Thanks to the geodesic condition of Bardos-Lebeau-Rauch (cf.~\cite{BLR92})
for the controllability of the wave equation,
the application of this method 
to the boundary controllability of the heat equation 
(cf.~\S\ref{sec:appli})
yields new geometric bounds on the cost of fast controls
(extending the results of \cite{LMheatcost} on internal controllability).
The companion paper~\cite{LMconscost} concerns 
the quite different case $\abs{\theta}=\pi/2$ 
(in particular there is no smoothing effect).

The relationship between first and second order controllable systems 
has been investigated in previous papers, 
always with $\theta=0$ and the additional initial data $\dz(0)=z_{1}$ in $X$
(i.e. considering trajectories of (\ref{eqsystcon}) 
in the state space $X\times X$).
In \cite{Fat67}, it is proved that 
the approximate controllability of (\ref{eqsystcon}) with $\dz(0)=z_{1}$
in some time implies the approximate controllability of (\ref{eqfsystcon})  
for any time (the control transmutation method yields an alternative proof), 
and proves the converse under some assumptions on the spectrum of $A$.
The converse is investigated further in~\cite{Tsu70} 
(in Hilbert spaces) and \cite{Tri78}. 
In a restricted setting,
the null controllability of (\ref{eqfsystcon}) was deduced from 
the exact controllability of (\ref{eqsystcon}) with $\dz(0)=z_{1}$
in \cite{Rus73} and \cite{Sei78}
by the indirect method of bi-orthogonal bases.

The study of the cost of fast controls 
was initiated by Seidman in~\cite{Sei84} 
with a result on the heat equation obtained by Russell's method
(cf. \cite{FCZ00}, \cite{LMheatcost} for improvements and other references).
Seidman also obtained results on the Schr{\"o}dinger equation 
by working directly on the corresponding 
window problem for series of complex exponentials 
(see \cite{LMschrocost} for improvements and references). 
With collaborators, he later treated 
the case of finite dimensional linear systems 
(cf. \cite{SY96}) 
and generalized the window problem to a larger class of complex exponentials 
(cf. \cite{SAI00}).
The control transmutation method generalizes 
upper bounds on the cost of fast controls 
from the one-dimensional setting 
(which reduces to a window problem) 
to the general setting which we specify in the next section.


\section{The setting}
\label{sec:setting}

We assume that $A$ is the generator of a 
strongly continuous cosine operator function $\Cos$
(i.e. the second-order Cauchy problem for $\ddz(t)=A(t)z$ 
is well posed and $\Cos$ is its propagator).
For a textbook presentation of cosine operator functions, 
we refer to chap.~2 of \cite{Fat85} or \S3.14 of \cite{ABHN01}. 
The associated sine operator function is $\Sin(t)=\int_{0}^{t}\Cos(s)ds$
(with the usual Bochner integral).
$\Cos$ and $\Sin$ are 
strongly continuous functions on $\R$ of bounded operators on $X$.
Moreover $A$ generates a holomorphic semigroup $T$ of angle $\pi/2$ 
(cf. th.~3.14.17 of \cite{ABHN01}).
In particular $S(t)=T(e^{i\theta}t)$ defines 
a strongly continuous semigroup $(S(t))_{t\in\Rp}$ of bounded operators on $X$.
In this setting, for any source term $f\in L^{1}\loc(\R,X)$, 
for any initial data $x_{0}$, $z_{0}$ and $z_{1}$ in $X$, 
the inhomogeneous first and second order Cauchy problems
\begin{align}
\label{eqfCpb}
&\dx(t)=e^{i\theta}A(t)x+f(t) & (t\in \Rp ), \quad  & x(0)=z_{0},\\
\label{eqCpb}
&\ddz(t)=A(t)z+f(t) & (t\in \R ), \quad  & z(0)=z_{0},\ z(0)=z_{1}, 
\end{align}
have unique mild solutions 
$x\in C^{0}(\Rp,X)$ and $z\in C^{0}(\R,X)$ 
defined by: 
\begin{gather*}
x(t)=S(t)x_{0} + \int_{0}^{t}\! S(t-s)f(s)dt ,\
z(t)=\Cos(t)z_{0} + \Sin(t)z_{1} + \int_{0}^{t}\! \Sin(t-s)f(s)dt .
\end{gather*}

\begin{remark}
When $A$ is a negative self-adjoint unbounded operator on a Hilbert space,
$T$, $\Cos$ and $\Sin$ are simply defined by the functional calculus as 
$T(t)=\exp(tA)$, $\Cos(t)=\cos(t\smA)$ and $\Sin(t)=(\smA)^{-1}\sin(t\smA)$.
\end{remark}

Following \cite{Weicon89}, 
we now make natural assumptions on $B$ 
for any initial data in the state space $X$
to define a unique continuous trajectory   
of each system (\ref{eqfsystcon}) and (\ref{eqsystcon}).
Let $X_{-1}$ be the completion of $X$ with respect to the norm 
$\norm{x}_{-1}=\norm{(A-\beta)^{-1}x}$
for some $\beta\in\C$ outside the the spectrum of $A$.
$X_{-1}$ is also the dual of the space $X_{1}$ defined 
as the domain $\D(A)$ with the norm $\norm{x}_{1}=\norm{(A-\beta)x}$.
We assume that $B\in\L(U,X_{-1})$ is 
{\em an admissible unbounded control operator} in the following sense:
\begin{gather}
\label{eqadmcon}
\forall t>0, \forall u\in L^{2}([0,t];U), \
\int_{0}^{t}S(s)Bu(s)ds \in X
\text{ and }
\int_{0}^{t}\Sin(s)Bu(s)ds \in X
.
\end{gather} 
In this setting, for any $x_{0}$ and $z_{0}$ in $X$,
for any $u$ and $v$ in $L^{2}\loc(\Rp;U)$,
the unique solutions $x$ and $z$ in $C^{0}(\R;X)$  
of (\ref{eqfsystcon}) and (\ref{eqsystcon}) respectively 
are defined by:
\begin{gather}\label{eqintcon}
x(t)=S(t)x_{0} + \int_{0}^{t}\! S(t-s)Bu(s)dt 
, \ 
z(t)=\Cos(t)z_{0}+\int_{0}^{t}\! \Sin(t-s)Bv(s)ds
. \end{gather}
The natural notions of controllability cost for the linear systems
(\ref{eqfsystcon}) and  (\ref{eqsystcon}) are:
\begin{definition}
\label{def:con}
The system (\ref{eqfsystcon}) 
is {\em null-controllable in time $T$} if 
for all $x_{0}$ in $X$, 
there is a $u$ in $L^{2}(\R;U)$ such that 
$u(t)=0$ for $t\notin [0,T]$ and $x(T)=0$.
The {\em controllability cost} for (\ref{eqfsystcon}) in time $T$ is
the smallest positive constant $\kappa_{1,T}$ 
in the following inequality for all such $\phi_{0}$ and $u$:
$\int_{0}^{T}\norm{u(t)}^{2}dt\leq \kappa_{1,T}\norm{x_{0}}^{2}$.

The system (\ref{eqsystcon}) 
is {\em null-controllable in time $T$} if 
for all $z_{0}$ in $X$,
there is a $v$ in $L^{2}(\R;U)$ such that 
$v(t)=0$ for $t\notin [0,T]$ and $z(T)=\dz(T)=0$.
The {\em controllability cost} for (\ref{eqsystcon}) in time $T$ is
the smallest positive constant $\kappa_{2,T}$ 
in the following inequality for all such 
$z_{0}$ and $v$:
$\int_{0}^{T}\norm{v(t)}^{2}dt\leq \kappa_{2,T}\norm{z_{0}}^{2}$.
\end{definition}

\begin{remark}
Equivalently, 
for all $x_{T}$ in $S(T)X$, 
there is a $u$ in $L^{2}(0,T;U)$ such that 
$x(0)=0$ and $x(T)=x_{T}$,
and, for all $z_{0}$ and $z_{T}$ in $X$,
there is a $v$ in $L^{2}(0,T;U)$ such that 
$\dz(0)=\dz(T)=0$, $z(0)=z_{0}$ and $z(T)=z_{T}$.
\end{remark}


\section{The results and the method}
\label{sec:res}

Our estimate of the cost of fast controls for (\ref{eqfsystcon}) builds,
through the control transmutation method, 
on the same estimate for a simple system of type (\ref{eqfsystcon}),
i.e. on a segment $[0,L]$
with Dirichlet ($N=0$) or Neumann ($N=1$) condition at the left end
controlled at the right end through a Dirichlet condition: 
\begin{gather}
\label{eq1dp}
\partial_{t}\phi =e^{i\theta}\partial_{s}^{2} \phi
\ \text{on}\ \left]0,T\right[\times \left]0,L\right[, \ 
\partial_{s}^{N}\phi\res{s=0} = 0, \ 
\phi\res{s=L}=
u,\
\phi\res{t=0} = \phi_{0}
. \end{gather}
With the notations of \S\ref{sec:setting}, $x=\phi$,
$A=\partial_{s}^{2}$ on $X=L^{2}(0,L)$
with $\D(A)=\{f\in H^{2}(0,L) \,|\, \partial_{s}^{N}f(0)=f(L)=0\}$,
$\norm{\cdot}_{1}$ with $\beta=0$ is 
the homogeneous Sobolev $\dot{H}^{2}(0,L)$ norm,
and $B$ on $U=\C$ is the dual of $C\in\L(X_{1};U)$  
defined by $Cf=\partial_{s}f(L)$.

It is well-known that the controllability of this system 
reduces by spectral analysis to classical results 
on nonharmonic Fourier series.
The following upper bound for the cost of fast controls,
proved in \S\ref{sec:fcs}, is an application of a 
refined result of Avdonin-Ivanov-Seidman in~\cite{SAI00}.
\begin{theorem}
\label{th:1d}
There are positive constants $\alpha$ and $\gamma$
such that, for all $N\in\{0,1\}$,
$L>0$, $T\in \left]0,\inf(\pi,L)^{2}\right]$,
the controllability cost $\kappa_{L,T}$ of the system (\ref{eq1dp})
satisfies: $\kappa_{L,T}\leq \gamma \exp(\alpha L^{2}/T)$.
\end{theorem}

This theorem leads to a definition
of the optimal fast control cost rate for (\ref{eq1dp}):
\begin{definition}
\label{defin:alpha}
The rate $\alpha_{*}$ is the smallest positive constant such that 
for all $\alpha > \alpha_{*}$
there exists $\gamma>0$ satisfying the property stated in theorem~\ref{th:1d}.
\end{definition}
\begin{remark}
Computing $\alpha_{*}$ is an interesting open problem and its solution 
does not have to rely on the analysis of series of complex exponentials.
The best estimate so far is 
$\alpha_{*}\in [ 1/2, 4(36/37)^{2}]$ for $\theta=0$ (cf. \cite{LMheatcost}).
\end{remark}

Our main result is a generalization of theorem~\ref{th:1d} to 
the first-order system (\ref{eqfsystcon}) 
under some condition on the second-order system (\ref{eqsystcon}):
\begin{theorem}
\label{th}
If the system (\ref{eqsystcon}) is null-controllable
for times greater than $L_{*}$,
then the system (\ref{eqfsystcon}) is null-controllable
in any time $T$.
Moreover, the controllability cost $\kappa_{1,T}$ of (\ref{eqfsystcon}) 
satisfies the following upper bound
(with $\alpha_{*}$ defined above):
\begin{gather}
\label{equb}
\limsup_{ T\to 0}  
T \ln \kappa_{1,T} \leq \alpha_{*} L_{*}^{2}
. \end{gather}
\end{theorem}
\begin{remark}
\label{rem:th}
The upper bound (\ref{equb}) means that the norm of an input function $u$
steering the system (\ref{eqfsystcon}) from an initial state $x_{0}$
to zero grows at most like $\gamma\norm{x_{0}}\exp(\alpha L^{2}/(2T))$
as the control time $T$ tends to zero
(for any $\alpha_{*}$ and some $\gamma>0$).
The falsity of the converse of the first statement in th.~\ref{th} 
is well-known, 
e.g. in the more specific setting of \S\ref{sec:appli}. 
\end{remark}
\begin{remark}
As observed in~\cite{Car93}, 
(\ref{equb}) yields a
logarithmic modulus of continuity for the minimal time function 
$\Tmin:X\to [0,+\infty)$ of (\ref{eqfsystcon}); i.e.
$\Tmin(x_{0})$, defined as the infimum of the times $T>0$ for which 
there is a $u$ in $L^{2}(\R;Y)$ such that 
$\int_{0}^{T}\norm{u(t)}^{2}dt\leq 1$,
$u(t)=0$ for $t\notin [0,T]$ and $x(T)=0$,
satisfies: for all $\alpha>\alpha_{*}$, 
there is a $c>0$ such that, 
for all $x_{0}$ and $x_{0}'$ in $X$
with $\norm{x_{0}-x_{0}'}$ small enough,
$\abs{\Tmin(x_{0}) - \Tmin(x_{0}')}\leq 
\alpha L_{*}^{2}/\ln(c/\norm{x_{0}-x_{0}'})$.
\end{remark}

It is well known that the semigroup $T$ 
can be expressed as an integral over the cosine operator function $\Cos$
(cf. the second proof of th.~3.14.17 in \cite{ABHN01}):
\begin{gather}
\label{eqWei}
\forall x\in X, \forall t\in \C \text{ s.t. } \abs{\arg t}<\pi/2 , \quad 
T(t)x=\int k(t,s) \Cos(s)x ds\ ,
\end{gather}
where $k(0,s)=\delta(s)$ and 
$k(t,s)=\exp(-s^{2}/(4t))/\sqrt{\pi t}$ for $\Re t>0$.
This equation has been referred to as 
the abstract Poisson or Weierstrass formula.
Starting with the observation that $k$ 
is the fundamental solution of the heat equation on the line,  
i.e. $k$ is the solution of $\partial_{t}k=\partial_{s}^{2}k$ 
with the Dirac measure at the origin as initial condition,
the {\em transmutation control method} consists in replacing 
the kernel $k$ in (\ref{eqWei})
by some  {\em fundamental controlled solution} on the segment $[-L,L]$ 
controlled at both ends (cf. (\ref{eqtrans})).
The one dimensional th.~\ref{th:1d} is used to construct this 
fundamental controlled solution in \S\ref{sec:fcs}
and the transmutation is performed in \S\ref{sec:transmut}.


\section{The fundamental controlled solution}
\label{sec:fcs}

This section begins with 
an outline of the standard application of \cite{SAI00}
to the proof of th.~\ref{th:1d}.
Following closely \S5 of \cite{LMheatcost},
the rest of the section outlines 
the construction of  a ``fundamental controlled solution'' $k$
in the following sense,
where $\mathcal{D}'(\mathcal{O})$ denotes the space of 
distributions on the open set $\mathcal{O}$
endowed with the weak topology, 
$\mathcal{M}(\mathcal{O})$ denotes the subspace of 
Radon measures on $\mathcal{O}$,
and $\delta$ denotes the Dirac measure at the origin:
\begin{definition}
\label{def:fundcontrol}
The distribution $k\in C^{0}([0,T]; \mathcal{M}(]-L,L[))$
is a fundamental controlled solution
for (\ref{eqk1}) at cost $(\gamma,\alpha)$ if 
\begin{align}
&\partial_{t}k = e^{-i\theta}\partial_{s}^{2}k  
\quad {\rm in }\ \mathcal{D}'(]0,T[\times ]-L,L[)\ , 
\label{eqk1} \\
&k\res{t=0}  =  \delta \quad {\rm and }\quad k\res{t=T}  =  0 \ ,
\label{eqk2} \\
&\|k\|_{L^{2}(]0,T[\times ]-L,L[)}^{2}
\leq \gamma e^{\alpha L^{2}/T } .
\label{eqk3}
\end{align}
\end{definition}

The operator $A$ defined at the beginning of \S\ref{sec:res}
is negative self-adjoint on the Hilbert space $L^{2}(0,L)$.
It has a sequence $\{\mu_{n}\}_{n\in \mathbb{N}^{*}}$
of negative decreasing eigenvalues 
and an orthonormal Hilbert basis $\{e_{n}\}_{n\in \mathbb{N}^{*}}$ 
in $L^{2}(0,L)$ of corresponding eigenfunctions.
Explicitly: 
$\sqrt{-\mu_{n}}=\left(n+\nu\right)\pi/L$
with $\nu=0$ for $N=0$ (Dirichlet) and $\nu=1/2$ for $N=1$ (Neumann).
First note that th.~\ref{th:1d} can be reduced to the case $L=\pi$ 
by the rescaling $(t,s)\mapsto (\sigma^{2}t,\sigma s)$ 
with $\sigma=L/\pi$.
In terms of the coordinates
$c=(c_{n})_{n\in \mathbb{N}^{*}}$ of $A^{N/2}f_{0}$ 
in the Hilbert basis $(e_{n})_{n\in \mathbb{N}^{*}}$
where $f_{0}$ is the initial state of the dual observability problem,
th.~\ref{th:1d} with $L=\pi$
reduces by duality to the following window problem:
$\exists \alpha>0$, $\exists \gamma>0$, $\forall T\in \left]0,\pi^{2}\right]$,
\begin{gather*}
\label{eq:windpb}
\forall c\in l^{2}(\N^{*}), \quad
\sum_{n\in \mathbb{N}^{*}}
\abs{c_{n}}^{2}
\leq  \gamma e^{\alpha\pi^{2}/T }
\int_{0}^{T}\abs{F(t)}^{2}dt
\text{ where } F(t)= \sum_{n=1}^{\infty} c_{n}e^{\exp(i\theta)\mu_{n} t}
. \end{gather*}
Since this results from th.~1 of \cite{SAI00}
with $\lambda_{n}\sim ie^{i\theta}n^{2}$  
as in \S5:2 of \cite{SAI00},
the proof of th.~\ref{th:1d} is completed.

Now we consider a system governed by the same equation as (\ref{eq1dp})
but on the twofold segment $[-L,L]$ controlled at both ends:
\begin{gather}
\label{eqfcsp}
\partial_{t}\phi - e^{-i\theta}\partial_{s}^{2} \phi=0
\quad {\rm in}\ ]0,T[\times ]-L,L[ ,\quad 
\phi\res{s=\pm L} =u_{\pm} ,\quad
\phi\res{t=0} = \phi_{0} 
, \end{gather}
with initial state $\phi_{0}\in  L^{2}(0,L)$, 
input functions $u_{-}$ and $u_{+}$ in $L^{2}(0,T)$.
As in proposition~5.1 of~\cite{LMheatcost},
applying th.~\ref{th:1d} with $N=0$ to the odd part of $\phi_{0}$ 
and with $N=1$ to the even part of $\phi_{0}$
proves that the controllability cost of (\ref{eqfcsp})
is not greater than the controllability cost of (\ref{eq1dp})
and therefore satisfies the same estimate stated in th.~\ref{th:1d}.
As in proposition~5.2 of~\cite{LMheatcost},
we may now combine successively 
the smoothing effect of (\ref{eqfcsp}) with no input (i.e. $u_{+}=u_{-}=0$) 
and this controllability cost estimate 
(plugged into the integral formula expressing $\phi$ in terms
of $\phi_{0}$ and $u_{\pm}=\phi\res{s=\pm L} $)
to obtain:
\begin{proposition}
\label{prop:fundcontrol}
For all $\alpha>\alpha_{*}$, there exists $\gamma>0$ such that 
for all $L>0$ and $T\in\, ]0,\inf(\pi/2,L)^{2}]$
there is a fundamental controlled solution
for (\ref{eqfcsp}) at cost $(\gamma,\alpha)$
(cf. def.~\ref{def:fundcontrol}).
\end{proposition}


\section{The transmutation of second-order controls into first-order controls}
\label{sec:transmut}

In this section we prove th.~\ref{th}.

Let $x_{0}\in X$ be an initial state for (\ref{eqfsystcon})
and let $L>L_{*}$.
Let $z\in C^{0}(\Rp;X)$ and $v\in L^{2}(\Rp;U)$ 
be the solution and input function 
obtained by applying the exact controllability of 
(\ref{eqsystcon}) in time $L$ to the initial state $z_{0}=x_{0}$.

We define $\zz\in C^{0}(\R;X)$
and $\vv\in L^{2}(\mathbb{R};U)$ 
as the extensions of $\zeta$ and $v$ 
by reflection with respect to $s=0$,
i.e. $\zz(s)=\zeta(s)=\zz(-s)$ and $\vv(s)=v(s)=\vv(-s)$ for $s\geq 0$.
They inherit from (\ref{eqintcon}):
\begin{gather}\label{eqintzz}
\zz(t)=\Cos(t)x_{0}+\int_{0}^{t}\! \Sin(t-s)B\vv(s)ds
. \end{gather}
Def.~\ref{def:con} of $\kappa_{2,L}$ 
implies the following cost estimate for $\vv$: 
\begin{gather} 
\label{eqvv}
\int \norm{\vv(s)}^{2}ds
= 2\int_{0}^{L} \norm{v(s)}^{2}ds
\leq 2\kappa_{2,L}\norm{x_{0}}^{2}
. \end{gather}
Since $\D(A)$ is dense in $X$,
there is a sequence $(x_{n})_{n\in N^{*}}$ in $\D(A)$
converging to $x_{0}$ in $X$.
Since $X_{1}$ is dense in $X_{-1}$,
there is a sequence $(f_{n})_{n\in \N^{*}}$ in $C^{1}(\R;X_{1})$
converging to $B\vv$ in $L^{2}(\mathbb{R};X_{-1})$.
For each $n\in \N^{*}$, let $\zz_{n}$ be defined in $C^{2}(\R;X)$ by:
\begin{gather*}
\zz_{n}(t)=\Cos(t)x_{n}+\int_{0}^{t}\! \Sin(t-s)f_{n}(s)ds
, \end{gather*}
which converges to $\zz(t)$ in $X$ for all $t$ due to (\ref{eqintzz}).
Since $\zz_{n}$ is a genuine solution of 
$\ddzz(t)=A\zz_{n}(t)+f_{n}(t)$ (cf. lem.~4.1 of \cite{Fat85}), 
we have for all $\varphi$ in $\D(A')$: 
\begin{gather*} 
s\mapsto \scd{\zz_{n}(s)}{\varphi}\in H^{2}(\R)
\quad\text{and}\quad  
\frac{d^{2}}{ds^{2}} \scd{\zz_{n}(s)}{\varphi}
= \scd{\zz_{n}(s)}{A'\varphi}+\scd{f_{n}(s)}{\varphi}
. \end{gather*}
Hence, 
$\displaystyle 
\scd{\zz_{n}(t)}{\varphi}= \scd{x_{n}}{\varphi}
+\int_{0}^{t}\! (t-s)\scd{\zz_{n}(s)}{A'\varphi}
+\int_{0}^{t}\! (t-s)\scd{f_{n}(s)}{\varphi}
$.
Passing to the limit, yields
$\displaystyle 
\scd{\zz(t)}{\varphi}= \scd{x_{0}}{\varphi}
+\int_{0}^{t}\! (t-s)\scd{\zz(s)}{A'\varphi}
+\int_{0}^{t}\! (t-s)\scd{B\vv(s)}{\varphi}
$. 
Therefore:
\begin{align} 
\label{eqzz}
& s\mapsto \scd{\zz(s)}{\varphi}\in H^{2}(\R)
\quad\text{and}\quad  
\frac{d^{2}}{ds^{2}} \scd{\zz(s)}{\varphi}
=\scd{\zz(s)}{A\varphi}+\scd{B\vv(s)}{\varphi} , \\
\label{eqzzL}
 &\scd{\zz(s)}{\varphi}=0 
\quad\text{and}\quad  
\frac{d}{ds} \scd{\zz(s)}{\varphi}=0
\quad \text{ for } \abs{s}=L .
\end{align}

Let $\alpha>\alpha_{*}$ and $T\in]0,\inf(1,L^{2})[$.  
Let $\gamma>0$ and 
$k\in C^{0}([0,T]; \mathcal{M}(]-L,L[))$ 
be the corresponding constant and fundamental controlled solution
given by proposition~\ref{prop:fundcontrol}.
We define $\kk\in C^{0}(\Rp; \mathcal{M}(\R))$ 
as the extension of $k$ by zero,
i.e. $\kk(t,s)=\bar{k}(t,s)$ on $[0,T]\times]-L,L[$
and $\kk$ is zero everywhere else.
It inherits from $k$ the following properties 
\begin{align}
&\partial_{t}\kk = e^{i\theta}\partial_{s}^{2} \kk
\quad {\rm in }\  \mathcal{D}'(]0,T[\times ]-L,L[)\ , 
\label{eqkk1} \\
&\kk\res{t=0}  =  \delta \quad {\rm and }\quad \kk\res{t=T}  =  0 \ ,
\label{eqkk2} \\
&\|\kk\|_{L^{2}(\Rp\times \R)}^{2}
\leq \gamma e^{\alpha L^{2}/T } .
\label{eqkk3}
\end{align}

The main idea of the proof is to use $\kk$
as a kernel to transmute $\zz$ and $\vv$
into a solution $x$ and a control $u$ 
for (\ref{eqfsystcon}). 
The transmutation formulas:
\begin{gather} 
\label{eqtrans}
x(t)=\int \kk(t,s)\zz(s)\, ds  
 \quad {\rm and } \quad 
\forall t>0, \ u(t)=\int \kk(t,s)\vv(s)\, ds 
, \end{gather}
define $x\in C^{0}(\Rp;X)$ and $u\in L^{2}(\Rp;U)$
since $\kk\in C^{0}(\Rp;\mathcal{M}(\R))\cap L^{2}(\Rp;L^{2}(\R))$, 
$\zz\in C^{0}(\mathbb{R};X)$ and $\vv\in L^{2}(\mathbb{R};U)$.
The property $(\ref{eqkk2})$ of $\kk$ implies 
$x(0)=x_{0}$ and $x(T)=0$.
Equations (\ref{eqzz}), (\ref{eqzzL}) and (\ref{eqkk1}) imply,
by integrating by parts:
\begin{gather*}
\forall \varphi\in \D(\A'), \ 
t\mapsto \scd{x(t)}{\varphi}\in H^{1}(\Rp)
\text{ and }
\frac{d}{dt} \scd{x(t)}{\varphi} = \scd{x(t)}{A'\varphi} + \scd{Bu(t)}{\varphi}
. \end{gather*}
This characterizes $x$ as the unique solution of (\ref{eqfsystcon})
in the weak sense (cf. \cite{Bal77}), 
which implies that $x$ and $u$ satisfy (\ref{eqintcon}).
Since $\int \norm{u(t)}^{2}dt\leq 
\int\!\! \int \abs{\kk(t,s)}^{2}ds\, dt \int \norm{\vv(s)}^{2}ds$,
(\ref{eqkk3}) and (\ref{eqvv}) 
imply the following cost estimate which completes the proof of th.~\ref{th}: 
\begin{gather*} 
\int_{0}^{T} \norm{u(t)}^{2}dt
\leq 2\kappa_{2,L}\gamma e^{\alpha L^{2}/T }  \norm{x_{0}}^{2}
. \end{gather*}


\section{Geometric bounds on the cost of fast boundary controls 
for the heat equation}
\label{sec:appli}

When the second-order equation (\ref{eqsystcon}) 
has a finite propagation speed
and is controllable, 
the control transmutation method yields 
geometric upper bounds on the cost of fast controls 
for the first-order equation (\ref{eqfsystcon}).
From this point of view, this method is an adaptation 
of the kernel estimates method of Cheeger-Gromov-Taylor in \cite{CGT82}.
This was illustrated in \cite{LMheatcost} and \cite{LMschrocost}
on the internal controllability of heat and Schr{\"o}dinger equations
on Riemannian manifolds
which have the wave equation as corresponding second-order equation.
Some similar lower bounds are proved in these papers
(without assuming the controllability of the wave equation)
which imply that 
the upper bounds are optimal with respect to time dependence.
In this section, we illustrate the control transmutation method 
on the analogous boundary control problems for the heat equation.

Let $(M,g)$ be a smooth connected compact 
$n$-dimensional Riemannian manifold with metric $g$ 
and smooth boundary $\partial M$. 
When $\partial M\neq\emptyset$, $M$ denotes the interior 
and $\adh{M}=M\cup\partial{M}$.
Let $\Delta$ denote the (negative) Laplacian on $(M,g)$
and $\partial_{\nu}$ denote the exterior Neumann vector field on $\partial M$.
The characteristic function of a set $S$ is denoted by $\chi_{S}$.

Let $X=L^{2}(M)$.
Let $A$ be defined by $Af=\Delta f$ 
on $\D(A)=H^{2}(M)\cap H^{1}_{0}(M)$.
Let $C$ be defined from $\D(A)$ to $U=L^{2}(\partial M)$
by $Cf=\partial_{\nu} f\res{\Gamma}$
where $\Gamma$ is an open subset of $\partial M$,
and let $B$ be the dual of $C$.
With this setting, 
(\ref{eqfsystcon}) with $\theta=0$ and (\ref{eqsystcon})
are the heat and wave equations 
controlled by the Dirichlet boundary condition on $\Gamma$.
In particular (\ref{eqsystcon}) writes:
\begin{gather}\label{eqwavecon}
\begin{split}
&\partial_{t}^{2} z-\Delta z =0\ \text{ on }\ \R_{t}\times M, 
\quad 
z=\chi_{\Gamma} v\ \text{ on }\ \R_{t}\times \partial M, \\
&z(0)=z_{0}\in L^{2}(M),\ \dza(0)=0,\  
v\in L^{2}\loc(\R;L^{2}(\partial M)) 
, \end{split}
\end{gather}
It is well known that $B$ is an admissible observation operator 
(cf. cor.~3.9 in \cite{BLR92}).
To ensure existence of a null-control for the wave equation
we use the geometrical optics condition of Bardos-Lebeau-Rauch
(specifically example~1 after cor.~4.10 in \cite{BLR92}):
\begin{gather}\label{eqGC}
\text{\parbox{\textequation}{
There is a positive constant  $L_{\Gamma}$ 
such that every generalized geodesic of length greater than $L_{\Gamma}$ 
passes through $\Gamma$ at a non-diffractive point.
}}
\end{gather} 
Generalized geodesics are the rays of geometrical optics 
(we refer to \cite{Mil02} for a presentation of this condition 
with a discussion of its significance).
We make the additional assumption
that they can be uniquely continued at the boundary $\partial M$.
As in~\cite{BLR92}, to ensure this, we may assume either that 
$\partial M$ has no contacts of infinite order with its tangents
(e.g. $\partial M=\emptyset$),
or that $g$ and $\partial M$ are real analytic.
For instance, we recall that (\ref{eqGC}) holds when 
$\Gamma$ contains a closed hemisphere 
of a Euclidean ball $M$ of diameter $L_{\Gamma}/2$,
or when $\Gamma=\partial M$ and $M$ is a strictly convex bounded Euclidean set
which does not contain any segment of length $L_{\Gamma}$. 
\begin{theorem}[\cite{BLR92}]
\label{th:BLR}
If (\ref{eqGC}) holds then the wave equation (\ref{eqwavecon}) 
is null-controllable in any time greater than $L_{\Gamma}$.
\end{theorem}

Thanks to this theorem, th.~\ref{th} implies:
\begin{theorem}
\label{th:geomub}
If (\ref{eqGC}) holds then the equation:
\begin{gather*}
\begin{split}
&\partial_{t}x-e^{i\theta}\Delta x =0\ \text{ on }\ \R_{t}\times M, 
\quad 
x=\chi_{\Gamma} u\ \text{ on }\ \R_{t}\times \partial M, \\
&x(0)=x_{0}\in H^{-1}(M),\ u\in L^{2}\loc(\R;L^{2}(\partial M)) 
, \end{split}
\end{gather*}
is null-controllable in any time $T$.
Moreover, 
the controllability cost $\kappa_{1,T}$ 
(cf. def.~\ref{def:con})
satisfies (with $\alpha_{*}$ as in def.~\ref{defin:alpha}):
$\displaystyle
\limsup_{ T\to 0}  
T \ln \kappa_{1,T} \leq \alpha_{*} L_{\Gamma}^{2}
$.
\end{theorem}




\end{document}